\documentclass{amsart}
\usepackage{amsmath,amsfonts}
\usepackage{amscd,amsthm}

\theoremstyle{plain} 

\newtheorem{theorem}{Theorem}[section]
\newtheorem{corollary}[theorem]{Corollary}
\newtheorem{lemma}[theorem]{Lemma}
\newtheorem{proposition}[theorem]{Proposition}

\newtheorem{remark}[theorem]{Remark}

\newcommand{\R}{\mathop{\mathbb{R}}}
\newcommand{\Q}{\mathop{\mathbb{Q}}}
\newcommand{\Z}{\mathop{\mathbb{Z}}}
\newcommand{\N}{\mathop{\mathbb{N}}}
\newcommand{\C}{\mathop{\mathbb{C}}}

\newcommand{\diag}[1]{\mathop{\it{diag}(#1)}}

\numberwithin{equation}{section}

\begin{document}

\title[Nekhoroshev--like estimate for non--linearizable analytic germs.]{Exponentially long time stability for non--linearizable analytic germs of $(\C^n,0)$.}

\author{Timoteo Carletti}

\date{\today}

\address[Timoteo Carletti]{Dipartimento di Matematica "U. Dini", viale Morgagni 67/A, 50134 Firenze, Italy}

\email[Timoteo Carletti]{carletti@math.unifi.it}


%
%

\begin{abstract}
We study the Siegel--Schr\"oder center problem on the linearization of analytic germs of 
diffeomorphisms in several complex variables, in the Gevrey--$s$, $s>0$ category. We 
introduce a new arithmetical condition of Bruno type on the linear part of the given germ, 
which ensures the existence of a Gevrey--$s$ formal linearization. We use this fact
to prove the effective stability, i.e. stability for finite but long time, of neighborhoods of the origin for the analytic germ.
\end{abstract}

\maketitle

\section{Introduction}

In this paper we consider the {\em Siegel--Schr\"oder center problem}~\cite{Herman,CarlettiMarmi,Carletti}
 in some class of ultradifferentiable germs of $(\C^n,0)$, $n \geq 1$; let us consider two 
classes of formal power series $\mathcal{A}_1 \subset \mathcal{A}_2 \subset \mathbb{C}^n \left[ \left[ z_1,\dots ,z_n \right] \right]$,
 closed w.r.t. to derivation and composition, let $F\in \mathcal{A}_1$ and 
call $DF(0)=A\in GL(n,\C)$, 
we say that $F$ is {\em linearizable in $\mathcal{A}_2$} if there exists $H\in \mathcal{A}_2$, normalized with
 $DH(0)=\mathbb{I}$, which solves~\footnote{Here $F\circ H$ means the composition of $F$ and
 $H$; in the following we will denote the composition of $F$ $n$-times with itself,
 by $F^n$ instead of $F^{\circ n}$.}:
\begin{equation}
  \label{eq:linearization}
  F \circ H (z)= H \circ R_{A}(z) \, ,
\end{equation}
where $R_{A}(z)=Az$. In the following we will assume $A$ to be diagonal with eigenvalues of unit modulus $\lambda_1,\dots,\lambda_n$, thus $A=\diag{\lambda_1,\dots,\lambda_n}$.

\indent
If both $\mathcal{A}_1$ and $\mathcal{A}_2$ coincide with the ring of formal power series then we have formal
 linearization if and only if $A$ is {\em non--resonant}, namely for all
 $\alpha \in \N^n$ such that $|\alpha|=\sum_{1\leq i \leq n}\alpha_i\geq 2$,
 and for all $j\in \{ 1,\dots,n \}$ then $\lambda^{\alpha}-\lambda_j \neq 0$ (where we 
used the standard notation $\lambda^{\alpha}=\lambda_1^{\alpha_1} \dots \lambda_n^{\alpha_n}$).

When $F$ is a germ of analytic diffeomorphisms defined in a neighborhood of the origin and we 
want to solve~\eqref{eq:linearization} in the same class of analytic germs, we have
to consider several cases. If $A$ is the {\em Poincar\'e domain}, namely $\sup_{1\leq j\leq n} |\lambda_j| <1$
 or $\sup_{1\leq j\leq n} |\lambda^{-1}_j| <1$, then Koenigs~\cite{Koenigs} and
 Poincar\'e~\cite{Poincare1} proved that every analytic germ $F \in Dif\!f(\C^n,0)$ 
such that $F(0)=0$ and $DF(0)=A$, is analytically linearizable. When $A$ is not 
in the Poincar\'e domain, we say that it is in the {\em Siegel domain}; the question is 
harder and some additional arithmetical conditions on $(\lambda_j)_j$ are needed
 (see~\cite{Herman} \S 17 page 158).

Let $p\in \N$, $p\geq 2$ and let us define for non--resonant $\lambda_1 ,\dots , \lambda_n$:
\begin{equation}
  \label{eq:smalldiv}
  \Omega(p)= \min_{1\leq j\leq n} \inf_{\substack{\alpha \in \N^n \\  0<|\alpha|<p}}|\lambda^{\alpha}-\lambda_j| \, ;
\end{equation}
we say that $A$ verifies a {\em Diophantine condition} of type $(\gamma,\tau)$ if 
there exist $\gamma>0$ and $\tau\geq 0$ such that for all $\beta\in \N^n\setminus \{ 0 \}$
 we have $\Omega(|\beta|)\geq \gamma |\beta|^{-\tau}$. Siegel~\cite{Siegel} in 1942 
 for the $n=1$ case and then Sternberg~\cite{Sternberg} and Gray~\cite{Gray} in the general
 case proved that if $A$ verifies a Diophantine condition then the linearization problem has
 an analytic solution. Bruno~\cite{Bruno} weakened the arithmetical condition by 
asking the convergence of the series $\sum_{k} \frac{\log \Omega^{-1}(2^{k+1})}{2^k}$. 
We remark that in the one dimensional case the Bruno condition~\footnote{In this case 
let $\omega \in (0,1)\setminus \Q$ such that $\lambda = e^{2\pi i \omega}$ and let $(q_n)_n$
be the denominators of the convergents~\cite{HardyWright} to $\omega$, then the Bruno
condition is equivalent to the convergence of the series 
$\sum_{k\geq 0}\frac{\log q_{k+1}}{q_k}$.} is optimal, as proved by
Yoccoz~\cite{Yoccoz}.

In~\cite{CarlettiMarmi} authors studied the Siegel--Schr\"oder center problem in the case of
general algebras of ultradifferentiable germs of $(\C,0)$, including the Gevrey case. In~\cite{Carletti} the multidimensional case is considered: 
if $\mathcal{A}_1=\mathcal{A}_2$ and $A$ verifies a Bruno 
condition, then every $F\in \mathcal{A}_1$ with $F(0)=0$ and $DF(0)=A$ is linearizable in $\mathcal{A}_2$, whereas if $\mathcal{A}_1$ is properly contained in $\mathcal{A}_2$ new conditions weaker than Bruno are sufficient to ensure linearizability in $\mathcal{A}_2$.

\indent
In this paper we consider in detail the case where $\mathcal{A}_1$ is the ring of germs of analytic 
diffeomorphisms at the origin of $n\geq 1$ complex variables, and $\mathcal{A}_2$ is the algebra of
 {\em Gevrey}--$s$, $s>0$, formal power series: the {\em Gevrey--$s$ linearization of 
analytic germs}.

 Let $\Hat F=\sum f_{\alpha} z^{\alpha}$, $(f_{\alpha})_{\alpha \in \N^n} \subset \C^n$ be a 
formal power series, then we say that it is {\em Gevrey--$s$}~\cite{Balser,Ramis1}, $s>0$, if there exist two
 positive constants $C_1,C_2$ such that:
\begin{equation}
  \label{eq:gevreydefvect}
  |f_{\alpha}| \leq C_1 C_2^{-s|\alpha|} |\alpha|!^s  \quad \forall \alpha \in \mathbb{N}^n \, .
\end{equation}
We denote the class of all formal vector valued power series Gevrey--$s$ by $\mathcal{C}_s$. It is closed w.r.t. derivation and composition.

\indent
In the Gevrey--$s$ case the arithmetical condition introduced in~\cite{CarlettiMarmi,Carletti} will be called {\em Bruno}--$s$ condition, $s>0$: for short $A\in \mathcal{B}_s$ if there exists
 an increasing  sequence of positive integer $(p_k)_k$ such that:
\begin{equation}
  \label{eq:brunosndim}
  \limsup_{|\alpha|\rightarrow +\infty}\left( 2\sum_{m=0}^{\kappa(\alpha)} \frac{\log \Omega^{-1}(p_{m+1})}{p_m}-s\log |\alpha|\right)< +\infty \, ,
\end{equation}
where $\kappa(\alpha)$ is defined by $p_{\kappa(\alpha)}\leq |\alpha| < p_{\kappa(\alpha)+1}$. 
When $n=1$ the Bruno--$s$ condition
 can be slightly weakened (see~\cite{CarlettiMarmi}); let $\omega \in (0,1)\setminus \Q$ and 
$\lambda=e^{2\pi i \omega}$, then the Bruno--$s$, $s>0$, condition reads:
\begin{equation}
\label{eq:brunos1dim}
\limsup_{n\rightarrow +\infty} \left( \sum_{j=0}^{k(n)}\frac{\log q_{j+1}}{q_{j}} - s\log n \right) <+\infty \, ,
\end{equation}
where $k(n)$ is defined by $q_{k(n)}\leq n < q_{k(n)+1}$. We remark that
 in both cases the new conditions are weaker than Bruno condition, which
is recovered when $s=0$. When $n=1$ we prove that the set $\mathcal{B}_s$ is $PSL(2,\Z)$--invariant (see remark~\ref{rem:invariance}). The main result of~\cite{Carletti} in the case of Gevrey--$s$ classes reads:
\begin{theorem}[Gevrey--$s$ linearization]
\label{thm:gevreylin}
Let $\lambda_1,\dots,\lambda_n$ be complex numbers of unit modulus and $A=\diag{\lambda_1,\dots,\lambda_n}$;
let $D_1 = \{ z \in \C^n : |z_i|<1 \, , 1\leq i\leq n \}$ be the isotropic polydisk 
of radius $1$ and let $F:D_1\rightarrow \C^n$ be an analytic function,
such that $F(z)=Az+f(z)$, with $f(0)=Df(0)=0$. If $A$ is non--resonant and verifies
a Bruno--$s$, $s>0$, condition~\eqref{eq:brunosndim} (or~\eqref{eq:brunos1dim} when $n=1$),
then there exists a formal Gevrey--$s$ linearization $\Hat{H}$ which 
solves~\eqref{eq:linearization}.
\end{theorem}
The aim of this paper is to show that the Gevrey character of the formal 
linearization can give information concerning the dynamics of the analytic germ.
Let $F(z)=Az+f(z)$ be a germ of analytic diffeomorphism verifying the hypothesis of Theorem~\ref{thm:gevreylin}, assume moreover $F$ not to be analytically 
linearizable. We will show that even if there is not {\em Siegel disk}, where the 
dynamics of $F$ is conjugate to the dynamics of its linear part, we have 
an open neighborhood of the origin which ``behaves as a Siegel disk'' under the 
iterates of $F$ for finite but long time, which results exponentially long:
the  {\em effective stability}~\cite{GDFGS} of the fixed point. 

In the case of analytic linearization, $|H_i^{-1}(z)|$, $i=1, \dots, n$, (which is well defined
 sufficiently close to the origin because $H$ is tangent to the identity) 
is {\em constant along the orbits}, namely it is a {\em first integral}
and $|F^m(z_0)|$ is bounded for all $m$ and sufficiently 
small $|z_0|$.

We will prove that any non--zero $z_0$ belonging to a polydisk of sufficiently
 small radius $r>0$, can be iterate a number of times 
$K=\mathcal{O}(exp \{ r^{-1/s} \} )$, being $s>0$ the Gevrey
exponent of the formal linearization, and we can find an 
{\em almost first integral}: 
a function which varies by a quantity of order $r$ under $m\leq K$ iterations, 
which implies that $F^m(z_0)$ is well defined and bounded for $m\leq K$. 
More precisely we prove the following

\begin{theorem}
  \label{thm:maintheorem}
Let $\lambda_1,\dots,\lambda_n$ be complex numbers of unit modulus and $A=\diag{\lambda_1,\dots,\lambda_n}$;
let $F:D_1\rightarrow \C^n$ be an analytic and univalent function,
such that $F(z)=Az+f(z)$, with $f(0)=Df(0)=0$. If $A$ is non--resonant and verifies
a Bruno--$s$, $s>0$, condition~\eqref{eq:brunosndim} (or~\eqref{eq:brunos1dim} when $n=1$),
then for all sufficiently small $0< r_{**} <1$, there exist positive constants 
$A_{**},B_{**},C_{**}$ such that for all $0<|z_0|<r_{**}/2$, the $m$--th iterate of 
$z_0$ by 
$F$ is well defined and  verifies $|z_m|=|F^m(z_0)|\leq C_{**} r_{**}$, for all 
$m\leq K_* = \Big \lfloor A_{**}^{-1}\, exp 
\Big \{ B_{**}\left( r_{**}/|z_0| \right)^{1/s} \Big \}\Big \rfloor$.
\end{theorem}
The hypothesis on the domain for $F$ is a natural normalization condition 
being the whole problem invariant by homothety.

In section~\ref{sec:conclusions} we compare our stability result with the stronger results which can be proved using Yoccoz's renormalization method~\cite{PerezMarco2} in the case $n=1$. Moreover we discuss the relation between our Bruno--$s$ condition and the arithmetical condition of P\'erez--Marco~\cite{PerezMarco1,PerezMarco2} ensuring that in the non--linearizable case the fixed point is accumulate by periodic orbits.

\indent
{\it Acknwoledgements.} 
I am  grateful to D. Sauzin for 
a very stimulating discussion concerning Gevrey 
classes and asymptotic analysis.

\section{Proof of the main Theorem}
\label{sec:proofmainthm}

In this part we will prove our main result, Theorem~\ref{thm:maintheorem}. The
proof will be divided into three steps: first we use the Gevrey--$s$ 
character of the formal linearization $\Hat H$, given by 
Theorem~\ref{thm:gevreylin}, to find an approximate solution of the 
conjugacy equation~\eqref{eq:linearization} up
to a (exponentially) small correction (paragraph~\ref{ssec:firststep});
then we prove an iterative Lemma allowing us to control how the small 
error introduced in the solution propagates (paragraph~\ref{ssec:thirdstep}). 
Finally we collect all
the informations to conclude the proof (paragraph~\ref{ssec:endproof}).

\subsection{Determination of an approximate solution}
\label{ssec:firststep}

We apply Theorem~\ref{thm:gevreylin}: the formal power series solution $\Hat H$ belongs to $\mathcal{C}_s$, as well as its inverse $\Hat{H}^{-1}$ which solves (formally):
\begin{equation}
  \label{eq:forhminus1}
  \Hat H^{-1} \circ F(z) = R_{A}\circ \Hat H^{-1}(z) \, .
\end{equation}
Since ${\Hat H}^{-1}=\sum h_{\alpha} z^{\alpha}\in \mathcal{C}_s$, there exist positive 
constants $A_1$ and $B_1$ such that 
\begin{equation}
\label{eq:gevreyhm1}
|h_{\alpha}| \leq A_1 B_1^{-s|\alpha|} |\alpha|!^s \quad \forall \,|\alpha| \geq 1 \, .
\end{equation}
For any positive integer $N$ we consider the {\em vectorial polynomial}, sum of
homogeneous vector monomials of degree $1\leq l \leq N$,
defined by: $\mathcal{H}_N(z)=\sum_{l=1}^{N} \sum_{|\alpha|=l} h_{\alpha} z^{\alpha}$ 
and the {\em Remainder Function}:
\begin{equation}
\label{eq:remainderfunction}
\mathcal{R}_N(z)=\mathcal{H}_N \circ F(z) - R_{A}\circ \mathcal{H}_N(z) \, .
\end{equation}
The following Proposition collects some useful properties of the remainder function.

\begin{proposition}
Let $\mathcal{R}_N(z)$ be the remainder function defined in~\eqref{eq:remainderfunction} and
 let $\alpha \in \N^n$, then:
  \begin{enumerate}
  \item[1)] $\partial_z^{\alpha} \mathcal{R}_N(0)=0$ if $|\alpha| \leq N$.
  \item[2)] For all $0<r<1$ there exists a positive constants $A_2$ and $B_2$
 such that if $|\alpha| \geq  N+1$, then:
\begin{equation*}
\Big |\frac{1}{\alpha !}\partial_z^{\alpha} \mathcal{R}_N(0)\Big |\leq A_2 
r^{-|\alpha|}B_2^{-sN} N!^s\, .
\end{equation*}
  \item[3)] For all $0<r<1$ and $|z|<r/2$ there exist positive constants $A_3,B_3$ such that:
\begin{equation}
\label{eq:punto3}
 |\mathcal{R}_N(z)| \leq A_3 B_3^{-sN} N!^s\left(\frac{|z|}{r}\right)^{N+1} \, .
\end{equation}
\end{enumerate}
Where we used the compact notation $\frac{1}{\alpha !}\partial_z^{\alpha} =
\frac{1}{\alpha_1 ! \dots \alpha_n !}\frac{\partial^{|\alpha|}}{\partial_{z_1}^{\alpha_1}
\dots \partial_{z_n}^{\alpha_n}}$.
\end{proposition}

\proof 
Statement 1) is an immediate consequence of the definition of $\mathcal{R}_N$.

To prove 2) we observe that $\mathcal{R}_N(z)$ is an analytic function on $D_1$ then one gets by Cauchy's estimates for all $0<r<1$ and for all $|\alpha|\geq N+1$:
\begin{equation}
  \label{eq:cauchy}
\Big | \frac{1}{\alpha!}\partial_z^{\alpha}\mathcal{R}_N(0)\Big |
\leq \frac{1}{(2\pi)^n} \frac{1}{r^{|\alpha|+1}}\max_{|z|=r}|\mathcal{H}_N\circ F(z)| \, .
\end{equation} 
Recalling the Gevrey estimate~\eqref{eq:gevreyhm1} for $\mathcal{H}_N$ and 
the analyticity of $F$ we obtain:
\begin{equation}
  \label{eq:estim2}
\Big | \frac{1}{\alpha!}\partial_z^{\alpha}\mathcal{R}_N(0)\Big |
\leq A_2 B_2^{-sN} N!^s r^{-|\alpha|} \, ,  
\end{equation}
for some positive constants $A_2$ and $B_2$ depending on the previous 
constants, on the dimension $n$ and on $F$.

To prove 3) we write the Taylor series $\mathcal{R}_N(z)=\sum_{|\alpha|\geq N+1}\frac{1}{\alpha!}\partial_z^{\alpha}\mathcal{R}_N(0) z^{\alpha}$: 
the bound on derivatives~\eqref{eq:estim2} 
implies
the estimate~\eqref{eq:punto3} for all $|z|< r/2$ and for some 
positive constants $A_3$ and $B_3$. 
\endproof

The bound~\eqref{eq:punto3} on $\mathcal{R}_N(z)$ depends on the positive 
integer $N$, so we can determine the value of $N$ for which the right
hand side of~\eqref{eq:punto3} attains its minimum, that's Poincar\'e's idea
of {\em summation at the smallest term}.

\begin{lemma}[Summation at the smallest term]
 \label{lem:sumupsmallest}

Let $\mathcal{R}_N(z)$ defined as before
and let $0<r_*<1/2$ then there exist positive constants $A_4,B_4$ such that 
for all $0<|z|<r_*$ we have:
\begin{equation}
  \label{eq:boundsumup}
  |\mathcal{R}_{\bar{N}}(z)|\leq A_4 \, exp \Big \{ -B_4
\left(\frac{r_*}{|z|}\right)^{1/s} \Big \} \, ,
\end{equation}
where $\bar{N}=\lfloor B_4\left(r_*/|z|\right)^{1/s}\rfloor$ and 
$\lfloor x \rfloor$ denotes the integer part of $x\in \R$.
\end{lemma}

\proof
Let us fix $0<r_*<1/2$, then for $0<|z|<r_*$ by Stirling formula
we obtain:
\begin{equation}
\label{eq:lem1}
| \mathcal{R}_N (z) | \leq A_4 \left( N B_3^{-1} \left(|z|/r_*\right)^{1/s} \right)^{Ns} e^{-sN} \, ,
\end{equation}
for some positive constant $A_4$. The right hand side of~\eqref{eq:lem1} attains its 
minimum at $\bar{N} = B_3 \left(r_*/|z|\right)^{1/s}$, evaluating the value of this minimum 
we get~\eqref{eq:boundsumup} with $B_4=B_3$.
\endproof

\subsection{Control of the ``errors''}
\label{ssec:thirdstep}

Let us define $\mathcal{H}(z)=\mathcal{H}_{\bar{N}}(z)$ and $\mathcal{R}(z)=
\mathcal{R}_{\bar{N}}(z)$, being $\bar{N}$ the ``optimal value'' obtained in
Lemma~\ref{lem:sumupsmallest}. We remark that $\mathcal{H}(z)$ doesn't 
solve~\eqref{eq:forhminus1} but the ``error'', $\mathcal{R}(z)$, is very 
small: exponentially small. We will prove that for initial conditions
in a sufficiently small disk, one can iterate an exponentially large
number of times without leaving a disk, say, of double size.

\begin{lemma}[Iteration lemma]
\label{lem:iteration}
Let $a,b,\alpha$ and $R$ be positive real numbers. Let us consider the sequence of positive 
number $(\mu_j)_{j\geq 0}$ defined by:
\begin{equation*}
\mu_0=R \quad \text{and} \quad \mu_{j+1} = 
\mu_j + a \, exp \{ -b/\mu_j^{\alpha} \} \, .
\end{equation*}
Let $K=\lfloor Ra^{-1}\, exp \{ b/(2R)^{\alpha} \}\rfloor$, then $\mu_j \leq 2R$ for all $j\leq K$.
\end{lemma}
\proof
Let us prove by induction on $j$ that for all $0\leq j\leq K$ we have 
\begin{equation}
  \label{eq:lem11}
  \mu_{j}\leq R + ja \, exp \{ -b/(2R)^{\alpha} \} \, ,
\end{equation}
then the claim will follow from~\eqref{eq:lem11} and the definition of $K$, 
in fact for all $j\leq K$:
\begin{equation*}
\mu_j \leq R + ja \, exp \{ -b/(2R)^{\alpha} \}\leq R+Ra^{-1} \, 
exp \{ b/(2R)^{\alpha} \} a\, exp \{ -b/(2R)^{\alpha} \}\leq 2R \, .
\end{equation*}
The basis of induction is easily verified; assume~\eqref{eq:lem11} 
for all $j \leq K-1$, we will prove it for $j=K$.
By definition of $(\mu_j)_j$ and the induction hypothesis we have:
\begin{eqnarray*}
  \mu_{K} &=& \mu_{K-1} + a \, exp \{ -b/\mu_{K-1}^{\alpha} \} \leq \\ 
  &\leq& R+(K-1)a \, exp \{ -b/(2R)^{\alpha}\}+a \, exp \{ -b/\mu_{K-1}^{\alpha} \} \, ,
\end{eqnarray*}
we remark that from~\eqref{eq:lem11} with $j=K-1$, using  $K-1 < Ra^{-1} \, 
exp \{ b/(2R)^{\alpha}\}$, we get $\mu_{K-1} \leq 2R$ and $exp \{ -b/\mu_{K-1}^{\alpha} \}
\leq exp \{ -b/(2R)^{\alpha} \}$. Then we conclude:
\begin{equation*}
  \mu_{K}\leq R +(K-1) a \, exp \{ -b/(2R)^{\alpha} \}+ a\, exp \{ -b/(2R)^{\alpha} \} \, ,
\end{equation*}
which ends the induction.
\endproof

Let $r_*$ as in Lemma~\ref{lem:sumupsmallest}, define $\rho(z)=|\mathcal{H}(z)|$ for 
all $0<|z|<r_*$, then Lemma~\ref{lem:sumupsmallest} admits the following Corollary,
which allows us to control the function $\rho(z)$ on consecutive points of an orbit of $F(z)$.

\begin{corollary}
  \label{lem:itarate}
Let $0<r_*<1/2$, let $r_1$ be the radius of the maximal polydisk where $\mathcal{H}(z)$ is invertible and let $r_{**}=\min (r_*, r_1)$. Then there exist 
positive constants $A_*,B_*$ such that for all $0<|z|<r_{**}$
 we have:
\begin{equation}
\label{eq:iterate}
\Big | \rho(F(z)) - \rho(z) \Big | \leq A_* \, exp \Big \{ -B_*
\left( \frac{r_*}{\rho(z)} \right)^{1/s} \Big \} \, .
\end{equation}
\end{corollary}

\proof
By definition $\rho(F(z))=|\mathcal{H}\circ F(z)|$ and $\rho(z)=|R_{A}\circ \mathcal{H}(z)|$,
since $|\lambda_j|=1$ for $1\leq j\leq n$, and $A=\diag{\lambda_1,\dots,\lambda_n}$, therefore:
\begin{equation*}
\Big |\rho(F(z)) - \rho(z) \Big | \leq \Big |\mathcal{H}\circ F(z) - R_{A}\circ \mathcal{H}(z) \Big | 
= |\mathcal{R}(z) |\, ,
\end{equation*}
and from Lemma~\ref{lem:sumupsmallest} we get:
\begin{equation}
\label{cor:eq1}
\Big |\rho(F(z)) - \rho(z) \Big | \leq A_4 \, exp \Big \{ -B_4 \left( 
\frac{r_*}{|z|}\right)^{1/s} \Big \} \, .
\end{equation}
We want to express this condition in terms of $\rho(z)$ instead of $|z|$, to do this
we have to consider the distortion properties of $\mathcal{H}(z)$ and of its inverse.
Let $J(z)=\partial_z \mathcal{H}(z)$ be the Jacobian of $\mathcal{H}(z)$ and 
let $J_1=\max_{|z|\leq r_1} |J(z)|$, where $r_1$ has been defined previously. 
Let $0<|z|<r_1$ and let us call 
$z^{\prime}=\mathcal{H}(z)$, clearly $|z^{\prime}|\leq J_1 r_1=r_2$. 
Let us call $J_2=\max_{|z^{\prime}|<r_2}|\partial_z \mathcal{H}^{-1}(z)|$, 
then for any $0<|z^{\prime}|<r_2$ there exists only one $z$ such that
$z=\mathcal{H}^{-1}(z^{\prime})$, which satisfies 
$|z|\leq J_2 |z^{\prime}|=J_2 |\mathcal{H}(z)|$.

Let $r_{**}=\min (r_*,r_1)$ then from~\eqref{cor:eq1} for any $0<|z|<r_{**}$
we get:
\begin{equation*}
\Big |\rho(F(z)) - \rho(z) \Big | \leq A_* \, exp \Big \{ -B_* 
\left( \frac{r_*}{\rho(z)}\right)^{1/s} \Big \} \, ,
\end{equation*}
where $A_*=A_4$ and $B_*=B_4J_2^{-1/s}$.
\endproof

\subsection{End of the proof}
\label{ssec:endproof}

We are now able to conclude the proof of the main 
Theorem~\ref{thm:maintheorem}. Take any $0<|z_0|<r_{**}/2$ and let us define 
$\rho_0=|z_0|$, $\rho_m=\rho(F^m(z_0))$ for all positive integer $m$ for
which $F^m(z_0)$ is well defined, by Corollary~\ref{lem:itarate} we have
\begin{equation}
  \label{eq:firstest}
  \rho_m \leq \rho_{m-1} + A_* \, exp 
\Big \{ -B_* \left( r_*/\rho_{m-1}\right)^{1/s} \Big \} \, .
\end{equation}
Let us call $R=|z_0|$, $a=A_*$, $b=B_* r_*^{1/s}$ and $\alpha=1/s$ then we 
can apply Lemma~\ref{lem:iteration} being $\mu_m \geq \rho_m$, to conclude
that:
\begin{equation}
  \label{eq:secondest}
  \rho_m \leq r_{**} \quad \forall m\leq K_*=\Big \lfloor |z_0|A_*^{-1} \, 
exp \Big \{ B_* \left( \frac{r_*}{2|z_0|}\right)^{1/s} \Big \}\Big\rfloor \, .
\end{equation}
This implies that $\mathcal{H}(z_m)$ is well defined in this range of values of
$m$, it is not constant and it evolves only by  
$\Big | |\mathcal{H}(z_m)| - |\mathcal{H}(z_0)| \Big |\leq r_{**}$. Recalling that
$z_m=F^m(z_0)$  we also have $|F^m(z_0)|\leq J_2 r_{**}$ and 
$\Big ||z_m|-|z_0| \Big | \leq J_2 r_{**}$ for all $0<|z_0|<r_{**}/2$ and all $m\leq K_*$.

This conclude the proof by setting $A_{**}=2A_*r_{**}^{-1}$, $B_{**}=B_* 
\left( r_*/(2r_{**})\right)^{1/s}$ and $C_{**}=J_2$.

\section{One dimensional case}
\label{sec:conclusions}

In this paper we proved that any analytic germs of diffeomorphisms of 
$(\C^n,0)$ with diagonal, non--resonant linear part has an 
{\em effective stability} domain, i.e. stable up to finite but ``long times'',
 close to the fixed point, provided the linear part verifies a new 
arithmetical Bruno--like 
condition~\eqref{eq:brunosndim} depending on a parameter $s>0$.

\begin{remark}[Invariance of $\mathcal{B}_s$, $n=1$ under the action of $PSL(2,\Z)$]
\label{rem:invariance}
The continued fraction development~\cite{HardyWright,MMY} of an irrational number $\omega$ gives us the sequences: $(a_k)_{k\geq 0}$ and $(\omega_k)_{k\geq 0}$. Then we introduce $(\beta_{k})_{k\geq -1}$ defined by $\beta_{-1}=1$ and for all integer $k\geq 0$: $\beta_{k}=\prod_{j=0}^k \omega_k$, which verifies : $1/2<\beta_kq_{k+1}<1$, where $q_k$'s are the denominators of the continued fraction development of $\omega$. We can then prove that condition Bruno--$s$~\eqref{eq:brunos1dim} is equivalent to the following one:
\begin{equation}
  \label{eq:brunosbeta}
  \limsup_{k \rightarrow +\infty}\left( \sum_{j=0}^k \beta_{j-1} \log \omega_j^{-1} + s \log \beta_{k-1} \right) < +\infty \, .
\end{equation}

Let us consider the generators of $PSL(2,\Z)$: $T\omega = \omega +1$ and $S\omega = 1/\omega$. For any irrational $\omega$, $T$ acts trivially being $\beta_k(T\omega)=\beta_k(\omega)$ for all $k$, whereas for irrational $\omega\in (0,1)$ we have $\beta_k(\omega)=\omega \beta_{k-1}(S\omega)$. Then the invariance of condition~\eqref{eq:brunosbeta} under the action of $GL(2,\Z)$ is obtained verifying the invariance under the action of $S$.
\end{remark}

Let us consider a slightly stronger version of the Bruno--$s$ condition: 
$\omega \in (0,1)\setminus \Q$ belongs to $\Tilde{\mathcal{B}}_s$ if:
\begin{equation}
    \label{eq:newbruno1}
  \lim_{n\rightarrow +\infty}\left( \sum_{l=0}^{k} \frac{\log q_{l+1}}{q_l}-s \log q_k \right) < +\infty \, ,
\end{equation}
where $(q_n)_n$ are the convergents to $\omega$. Let us introduce two other arithmetical conditions; 
let us denote by $\mathcal{B}_s^{\prime}$
 the set of irrational numbers whose convergents verify: 
\begin{equation}
  \label{eq:brunoprimes}
 \lim_{k\rightarrow +\infty} \frac{\log q_{k+1}}{q_k\log q_k}= s \, .
\end{equation}
The second condition is as follows, let $(\gamma_m)_{m\geq 1}$ and $(s_m)_{m\geq 1}$ 
be two positive sequences of real numbers such that:
$\sum_1^{+\infty} \gamma_m =\gamma < +\infty$ and $\sum_1^{+\infty} s_m =\sigma < +\infty$,
then we define a condition $\mathcal{B}_{\gamma,\sigma}$ by: 
\begin{equation}
  \label{eq:brunogamsig}
\frac{\log q_{m+1}}{q_m}\leq s_m \log q_m +\gamma_m \quad \forall m\geq 1\, .
\end{equation}

\begin{proposition}
  \label{prop:differentbruno1s}
Let $\omega \in (0,1)\setminus \Q$ and let $s>0$ then we have the following inclusions:
\begin{enumerate}
\item[1)] let $\omega \in \Tilde{\mathcal{B}}_s$, if $\omega$ is not a Bruno number
then $\omega \in \mathcal{B}_s^{\prime}$, otherwise $\omega \in \mathcal{B}^{\prime}_0$.
\item[2)] Let $\sigma \leq s$ and $\omega \in \mathcal{B}_{\gamma,\sigma}$ then 
$\omega \in \Tilde{\mathcal{B}}_s$;
\end{enumerate}
\end{proposition}

\proof
To prove the first statement let us write the following identity:
\begin{equation}
  \label{eq:step1}
  \sum_{l=0}^k \frac{\log q_{l+1}}{q_l}-s\log q_k = C + \sum_{l=2}^k 
\left[ \frac{\log q_{l+1}}{q_l}-s\left( \log q_l - \log q_{l-1} \right) \right] \, ,
\end{equation}
where $C=(1-s)\log q_1+\frac{\log q_{2}}{q_1}$. By condition $\Tilde{\mathcal{B}}_s$, 
this series converges and then its generic term goes to zero, from which we get:
\begin{equation}
  \label{eq:step3}
  \lim_{k\rightarrow +\infty}\frac{\log q_{k+1}}{q_k\log q_k}=s\left( 1-\lim_{k\rightarrow +\infty} \frac{\log q_{k-1}}{\log q_k} \right)\, .
\end{equation}
Let us denote by $s^{\prime}$ be value of the right hand side of~\eqref{eq:step3}, then clearly $s^{\prime}\in \left[ 0, s\right]$. Let us suppose $s^{\prime}>0$, but then we have for all sufficiently large 
$k$:
\begin{equation*}
  \frac{C_1}{q_{k}}\leq \frac{\log q_k}{\log q_{k+1}} \leq \frac{C_2}{q_{k}} \, ,
\end{equation*}
for some positive constants $C_1,C_2$, from which we get 
$\frac{\log q_k}{\log q_{k+1}}\rightarrow 0$, and from~\eqref{eq:step3} we conclude that 
$s^{\prime}=s$.

If $s^{\prime}=0$, namely $\frac{\log q_k}{\log q_{k+1}}\rightarrow 1$,
then it is easy to check that $\omega$ is a Bruno number.

Let us prove the second statement. For any positive integer $k$, using the definition
of $\mathcal{B}_{\gamma,\sigma}$ we can write:
\begin{equation}
  \sum_{l=0}^k \frac{\log q_{l+1}}{q_l}-s \log q_k \leq \sum_{l=0}^k s_l \log q_l +\sum_{l=0}^k \gamma_l - s \log q_k \, ,
\label{eq:primaeq}
\end{equation}
for all $0\leq l \leq k$ we have $\log q_l \leq \log q_k$ then the right hand side of~\eqref {eq:primaeq} is bounded by: $- \log q_k \left( s-\sum_{l=0}^k s_l \right)+\sum_{l=0}^k \gamma_l$. By hypothesis $\sum_{l=0}^k s_l \leq s$, for all $k$, then using $-\log q_k \leq -\log q_1$ we obtain:
\begin{equation*}
  \sum_{l=0}^k \frac{\log q_{l+1}}{q_l}-s \log q_k \leq - \log q_1 \left( s-\sum_{l=0}^k s_l \right)+\sum_{l=0}^k \gamma_l \, ,
\end{equation*}
then passing to the limit on $k$ we have:
\begin{equation*}
  \lim_{k\rightarrow +\infty} \sum_{l=0}^k \frac{\log q_{l+1}}{q_l}-s \log q_k \leq - \log q_1 \left( s-\sigma \right) + \gamma <+\infty \, .
\end{equation*}
\endproof

\begin{remark}
These new arithmetical conditions are weaker than the Bruno one, for instance 
condition $\mathcal{B}^{\prime}_s$ is verified by numbers $\omega$ whose denominators 
$(q_k)_k$ satisfy a growth condition like $q_{k+1}\sim q_k!^s$. Condition $\mathcal{B}_{\gamma,\sigma}$ implies convergence of the series:
$\sum_{k\geq 0}\frac{\log q_{k+1}}{q_k \log q_k}$.
\end{remark}

Let us conclude recalling a stability result of P\'erez--Marco~\cite{PerezMarco1,PerezMarco2}
and compare it with our result. In~\cite{PerezMarco2} author proved (Theorem V.2.1 Annexe 2
 \S f ) using a {\em geometric renormalization} scheme "\`a la Yoccoz" valid in the one dimensional case, 
a stability result that can be stated as follows:
\begin{theorem}[P\'erez--Marco, Contr\^ole de la diffusion]
\label{thm:pm}
Let $\omega \in (0,1)\setminus \Q$ and let $(q_k)_k$ be the denominators of its convergents.
Let $F$ be an analytic and univalent function defined in the unit disk 
$\{ z\in \C: |z|<1 \}$ such that $F(z)=\lambda z + \mathcal{O}(|z|^2)$, where 
$\lambda=e^{2\pi i \omega}$. There exist two positive constants $C_1,C_2$ such that if:
\begin{equation}
\label{eq:pm1}
|z| \leq C_1 e^{-\sum_{j=0}^{k-1} \frac{\log q_{j+1}}{q_j}} \, ,
\end{equation}
then for all integer $0\leq m \leq q_k$ we have:
\begin{equation}
\label{eq:pm2}
|F^m(z)| \leq C_2 e^{-\sum_{j=0}^{k-1} \frac{\log q_{j+1}}{q_j}} \, .
\end{equation}
\end{theorem}

The meaning of the Theorem is clear: if we start inside a disk of radius 
$r=C_1 e^{-\sum_{j=0}^{k-1} \frac{\log q_{j+1}}{q_j}}$ then we can apply $F$, up to
$q_k$ times, without leaving a disk of radius $rC_2/C_1$. To compare this result with our
effective stability result we have to make explicit the relation w.r.t. $r$ and $q_k$, 
which give the time of "stability".
Using our Bruno--$s$ condition~\eqref{eq:newbruno1} we can say that 
$C\leq rq^s_{k-1} \leq C^{\prime}$ for some positive constants $C,C^{\prime}$.
But from~\eqref{eq:brunoprimes} we get $\log q_k \leq C_3 q_{k-1}\log q_{k-1}$ for some
positive constant $C_3$, namely there exist positive constants $C^{\prime}_3,C_4$ such that:
\begin{equation*}
q_k \leq exp \Big \{ \frac{C^{\prime}_3}{r^{1/s}} \log \frac{C_4}{r^{1/s}} \Big \} \, .
\end{equation*}
We can then restate Theorem~\ref{thm:pm} as follows: if $|z|\leq r$, then 
$|F^m(z)|\leq rC_2/C_1$ for all integer 
$0\leq m \leq exp \{ \frac{C^{\prime}_3}{r^{1/s}} \log \frac{C_4}{r^{1/s}} \}$, obtaining
a better estimate on the time of effective stability.

We end with a last remark related again to the work of P\'erez--Marco.
\begin{remark}
P\'erez--Marco proved in~\cite{PerezMarco1,PerezMarco2} that any non--analytically
linearizable analytic germ, univalent in the unit disk, whose multiplier at the 
fixed point, verifies the following arithmetical condition:
\begin{equation}
\label{eq:perezemarco}
\sum_{k\geq 0}\frac{\log \log q_{k+1}}{q_k}<+\infty \, ,
\end{equation}
has a sequence of periodic orbits accumulating the fixed point, whose periods, $(q_{n_k})_k$,
make the Bruno series diverging.

Our Bruno--$s$ condition implies~\eqref{eq:perezemarco}, in fact 
from~\eqref{eq:newbruno1} we get:
\begin{equation*}
\sum_{k= 0}^N\frac{\log \log q_{k+1}}{q_k} \leq \sum_{k= 0}^N \left(\frac{\log C_3}{q_k}+
\frac{\log q_k}{q_k}+\frac{\log\log q_k}{q_k}\right) \, ,
\end{equation*}
we can let $N$ grow and using standard number theory results concerning the convergents, 
we obtain the P\'erez--Marco condition.
Then we can suppose these periodic orbits accumulating the fixed point to ``produce the 
effective stability: preventing the orbits from a too fast escape'', a situation similar to the one
holding in the Nekhoroshev Theorem for Hamiltonian systems~\cite{Nekhoroshev}, where the 
resonant web confines the flow for exponentially long times. It would be very interesting to know whether a similar phenomenon takes place in higher dimension.

We conclude by pointing out that our method gives us a stability exponent depending
on the Gevrey exponent and independent of the dimension: the bigger is the exponent, longer is
the time interval of stability, we can always take $s$ small enough to have a very
long time of stability. 
\end{remark}

\end{document}